\documentclass[11pt,oneside]{article}
\usepackage{amssymb,amsmath}
\usepackage{url}
\usepackage[english]{babel}
\usepackage{pstricks,pst-node}
\renewcommand{\le}{\leqslant}
\renewcommand{\ge}{\geqslant}
\newcommand{\bad}{\mathbf{Bad}}
\newcommand{\mad}{\mathbf{Mad}}

\newcommand{\Ll}{{\rm L}}
\newcommand{\area}{\mathbf{area}}
\newcommand{\vol}{\mathbf{vol}}

\newcommand{\RR}{\mathbb{R}}
\newcommand{\ZZ}{\mathbb{Z}}
\newcommand{\QQ}{\mathbb{Q}}
\newcommand{\KK}{\mathbf{K}}

\newcommand{\NN}{\mathbb{N}}

\newcommand{\BB}{\mathbf{B}}
\newcommand{\JJ}{\mathcal{J}}
\newcommand{\II}{\mathcal{I}}
\newcommand{\LLL}{\mathcal{L}}
\newcommand{\DDD}{\mathcal{D}}
\newcommand{\PPP}{\mathcal{P}}

\newcommand{\KKK}{\mathbf{K}}
\newcommand{\RRR}{\mathbf{R}}
\newcommand{\SSS}{\mathcal{S}}

\newcommand{\vr}{\mathbf{r}}

\newtheorem{lemma}{Lemma}
\newtheorem{theorem}{Theorem}
\newtheorem{proposition}{Proposition}
\newtheorem{problem}{Problem}

\newtheorem{conjecture}[problem]{Conjecture}

\newtheorem{vsec}{Proposition}

\newtheorem{BVThm}{Theorem (BV4)}

\newtheorem{BVThmm}{Theorem (BV5)}

\newenvironment{proof}{\textsc{Proof.}}{\ \newline\hspace*{\fill}$\boxtimes$}

\begin{document}

\title{On multiplicatively badly approximable numbers}

\author{Dzmitry Badziahin \footnote{Research supported by EPSRC grant EP/E061613/1}
\\ {\small\sc York}}

%\date{\it Dedicated to Graham Everest}

\maketitle

\begin{abstract}
The Littlewood Conjecture states that  $\liminf_{q\to \infty}
q\cdot||q\alpha||\cdot||q\beta||=0$  for all $(\alpha,\beta) \in
\RR^2$. We show that with the additional factor of $\log q\cdot
\log\log q$ the statement is false. Indeed,  our main result implies
that the set of $(\alpha,\beta)$ for which $ \liminf_{q\to\infty}
q\cdot\log q\cdot\log\log q\cdot||q\alpha||\cdot||q\beta||>0$ is of
full dimension.

\end{abstract}

\section{Introduction}

The famous Littlewood conjecture (LC) states that for any pair of
real numbers $(\alpha,\beta)$
\begin{equation}\label{little}
\liminf_{q\to \infty} q\cdot ||q\alpha||\cdot||q\beta||=0
\end{equation}
where $\|\cdot\|$ denotes the distance to the nearest integer.
Equivalently, the set
\begin{equation}\label{cond_lit}
\{(\alpha,\beta)\in \RR^2\;:\; \liminf_{q\to \infty}
q\cdot||q\alpha||\cdot||q\beta||>0\}
\end{equation}
is empty. This problem was conjectured in 1930's and it is still
open. For recent progress concerning this fundamental problem see
\cite{ELK, PVL} and references therein. It is easily seen that
\eqref{little} holds for all $\alpha\in\RR$ and $\beta\in \RR$
outside the set $\bad$ of badly approximable numbers defined as
follows
$$
\bad:=\{\alpha\in \RR\;:\; \liminf_{q\to \infty} q||q\alpha||>0\}.
$$

In attempt to understand what should be a proper analogue of badly
approximable points in multiplicative case several authors
investigated the following set (we will follow the notation
introduced in \cite{BV_mix}).
%
%In this paper we will investigate weakening the condition in
%\eqref{cond_lit} so as to guarantee that the corresponding set is
%nonempty. We will follow the notation introduced in \cite{BV_mix}.
For $\lambda\ge 0$ let
$$
\mad^{\lambda}:=\{(\alpha,\beta)\in \RR^2\;:\; \liminf_{q\to\infty}
\,(\log q )^\lambda \cdot q \cdot ||q\alpha||\cdot ||q\beta||>0\}.
$$
In other words, $\mad^\lambda$ is a modification of the set in
\eqref{cond_lit} such that the corresponding condition is weakened
by $(\log q)^\lambda$. More generally, given a function
$f\;:\;\NN\to \RR^+$, define the set
\begin{equation}\label{def_mp}
\mad(f):= \inf\{(\alpha,\beta)\in\RR^2\;:\; \liminf_{q\to
\infty}f(q)\cdot q \cdot ||q\alpha||\cdot||q\beta||>0\}.
\end{equation}

In \cite{BV_mix} the author and Velani conjectured that
\begin{conjecture}[BV]
\begin{eqnarray*}
\mad^\lambda=\emptyset\quad \text{for any }\lambda<1,\\
\dim(\mad^\lambda)=2\quad \text{for any }\lambda\ge 1\;\,
\end{eqnarray*}
\end{conjecture}
where $\dim(\cdot)$ denotes the Hausdorff dimension. If true this
conjecture implies that the proper multiplicative analogue of the
set $\bad$ is $\mad^1$. Note that LC is equivalent to the statement
that $\mad^0$ is empty. Therefore BV~conjecture implies LC.
Regarding the first part of BV~conjecture all that is known to date
is the remarkable result of Einsiedler, Katok and Lindenstrauss
\cite{ELK} which states that $\dim\mad^0=0$. On the other hand
according to the second part the best known result is due to Bugeaud
and Moschevitin \cite{BM}. It states that $\dim\mad^2=2$. So we have
a gap $ 0\le\lambda<2 $ where the behavior of $\mad^\lambda$ is
completely unknown.

In this paper we will address the second part of the BV~conjecture.
In particular, we will show that
$$
\dim\mad(f)=2\quad\text{if }\ f(q)=\log q\cdot \log\log q.
$$
It will straightforwardly imply that $\dim(\mad^\lambda)=2$ for any
$\lambda>1$.

 It is worth mentioning that the `mixed' analogue of this result
was achieved recently by author and Velani. It was proven that the
set
$$
\mad_\DDD(f):=\{\alpha\in\RR\; : \; \liminf_{q\to\infty} f(q)\cdot
q\cdot |q|_\DDD ||q\alpha||>0\}
$$
has full Hausdorff dimension. All the details can be found in
\cite{BV_mix}.

\subsection{Simultaneous and dual variants of $\mad$}

It is well known that Littlewood conjecture has an equivalent
formulation in terms of linear forms. In other words, \eqref{little}
is equivalent to the statement that
$$
\liminf_{|AB|\to\infty}|A|^*|B|^*\cdot||A\alpha-B\beta||>0
$$
where $|x|^*:=\max\{|x|,1\}$. However it is not known if
\eqref{def_mp} can be reformulated in the same manner. In other
words, define the sets
\begin{equation}\label{def_ml}
\mad_L(f):= \inf\{(\alpha,\beta)\in\RR^2\;:\; \liminf_{|AB|\to\infty
}f(|A|^*|B|^*) \cdot |A|^*|B|^*||A\alpha-B\beta||>0\}
\end{equation}
and
$$
\mad_L^\lambda:=\mad_L(\log^\lambda q).
$$
Then $\mad(f)$ and $\mad_L(f)$ are not necessarily the same. However
as it will be shown in the next sections these sets are closely
related to each other. For consistency in further discussion we will
use the notation $\mad_P^\lambda$ and $\mad_P(f)$ instead of
$\mad^\lambda$ and $\mad(f)$ respectively. It will reflect the fact
that in one case we deal with points and in another case we deal
with lines.

It appears that instead of investigating $\mad_P(f)$ and $\mad_L(f)$
independently it is easier to deal with them simultaneously. In
particular, we prove the following result:

%We call the pair $(\alpha,\beta)\in\RR^2$ multiplicatively badly
%approximable if $\omega (\alpha,\beta)\le 1$. Denote the set of
%multiplicatively badly approximable numbers by $\mad$. We believe
%that for all $(\alpha,\beta)\in \RR^2$, $\omega(\alpha,\beta)\ge 1$
%so in view of this the pairs $(\alpha,\beta)\in\mad$ are
%multiplicatively approximated by rationals in the worst possible
%way.
%
%In this paper we show that the set $\mad$ is nonempty and moreover
%it has full Hausdorff dimension. More exactly we prove the following
%result
\begin{theorem}\label{th_main1}
Let $f(q)=\log q\cdot \log\log q$. Then
$$
\dim(\mad_P(f)\cap\mad_L(f))=2.
$$
\end{theorem}

%Inparticular it significally improves the previous best known result
%of Moschevitin. And it is close to what we believe should be sharp
%result:
%\begin{conjecture}
%$$
%\mad_P(\log^\omega q)=\emptyset\quad \mbox{ if } \omega< 1
%$$
%and
%$$
%\dim(\mad_P(\log^\omega q))=2\quad\mbox{ if } \omega\ge 1.
%$$
%\end{conjecture}
%
%For convenience let's denote the set from the theorem by $\madl$.
%$$
%\madl:=\{(\alpha,\beta)\in \RR^2\;:\; \liminf_{q\to\infty} q\log
%q\log\log q\cdot||q\alpha||\cdot||q\beta||>0\}.
%$$
%One can see that $\madl\subset\mad$ therefore the statement of
%Theorem \eqref{th_main1} is even stronger than just $\dim\mad=2$. It
%also significally improves the the best previously known result of
%Moschevitin.

\subsection{Main result}

For convenience, we define the `modified logarithm' function
$\log^*\;:\; \RR\to\RR$ as follows
$$
\log^* x:=\left\{\begin{array}{ll}1&\mbox{for
}x<e;\\
\log x&\mbox{for }x\ge e.
\end{array}\right.
$$
%In the further discussion we will call by $f(q)$ the function
From now on
$$
f(q):=\log^* q\cdot\log^*\log q.
$$

The key to establishing Theorem \ref{th_main1} is to investigate the
intersection  of the sets $\mad_P(f)$ and $\mad_L(f)$ along fixed
vertical lines in the $(x,y)$-plane. With this in mind,  let $\Ll_x$
denote the line parallel to the $y$-axis passing through the point
$(x,0)$.

The following constitutes our main theorem.
\begin{theorem}\label{tsc}
For any $\theta \in \bad$
$$
\dim ( \mad_P(f)\cap\mad_L(f)\cap \Ll_\theta) = 1 \ .
$$
\end{theorem}

\noindent Since by Jarn\'{\i}k (1928) the Hausdorff dimension of
$\bad$ is one, Theorem~\ref{th_main1} can be easily derived from
Theorem~\ref{tsc} with the help of the following general result that
relates the dimension of a set to the dimensions of parallel
sections, enables us to establish the complementary lower bound
estimate  -- see \cite[pg. 99]{falc}.

\begin{vsec}
 Let $F$ be a subset of $\RR^2$ and let $E$ be
a subset of the $x$-axis. If \ $\dim (F \cap \Ll_x)  \ge t $ for all
$x \in E$, then $\dim  F  \ge  t + \dim E$.
\end{vsec}

\noindent Indeed, let $ F=\mad_P(f)\cap \mad_L(f) $ and $E=\bad$. In
view of $\dim(\bad)=1$ and Theorem~\ref{tsc}, one gets $\dim
\mad_P(f)\cap \mad_L(f) \ge 2$. Since $\mad_P(f)\cap \mad_L(f)
\subset \RR^2$, the upper bound statement for the dimension is
trivial. Therefore the main ingredient in establishing
Theorem~\ref{th_main1} is Theorem \ref{tsc}.

Regarding the proof of Theorem~\ref{tsc} we will use ideas similar
to those in \cite{BV_mix} which firstly appeared in joint work of
author, Pollington and Velani~\cite{BPV}. However the technical
details in this paper are substantially more complicated than those
in \cite{BV_mix}.

\section{Preliminaries}

Let $S$ be any subset of $\RR^2$. By $S_{\theta}$ we denote its
orthogonal projection onto the line $\Ll_\theta$. Let
$P(p,r,q):=(p/q,r/q)$ be a rational point where $(p,r,q)\in\ZZ^3,
\gcd(p,r,q)=1$. Denote by the height of $P$ the value
$$
H(P):=q^2|q\theta-p|\ge q^2||q\theta||.
$$
Denote by $\Delta(P,\delta)$ the following segment on $\Ll_\theta$:
$$
\Delta(P,\delta):=\{\theta\}\times\left(\frac{r}{q}-\frac{\delta}{H(P)},\frac{r}{q}+\frac{\delta}{H(P)}\right).
$$
So $|\Delta(P,\delta)|=2\delta H(P)^{-1}$.

Given a line with integer coeffitients
$$
L(A,B,C):=\{(x,y)\in \RR^2\;:\, Ax-By+C=0\},
$$
\begin{equation}\label{cond_abc}
(A,B,C)\in\ZZ^3,\ B\neq 0,\ \gcd(A,B,C)=1
\end{equation}
 denote by the height of $L$ the value
$$
H(L):=|A|^*B^2.
$$
Denote by $\Delta(L,\delta)$ the following segment on $\Ll_\theta$:
$$
\Delta(L,\delta):=\{\theta\}\times\left(
\frac{A\theta+C}{B}-\frac{\delta}{H(L)},\frac{A\theta+C}{B}+\frac{\delta}{H(L)}\right).
$$
So $|\Delta(L,\delta)|=2\delta H(L)^{-1}$.

Given constants $c>0$ and $Q>0$ define the auxiliary sets:
$$
\mad_P(f,c,Q):= \left\{(\alpha,\beta)\in \RR^2 \;:\; f(q)\cdot
q\cdot ||q\alpha||\cdot ||q\beta||>c\;\ \forall q\in\NN,\;
\ge Q\right\}
$$
and
$$
\mad_L(f,c,Q):= \inf\left\{(\alpha,\beta)\in\RR^2\;:\begin{array}{l}
f(|A|^*|B|^*)
\cdot |A|^*|B|^*||A\alpha-B\beta||>c,\\[1ex]
\forall (A,B)\in\ZZ^2,\ |A|^*B^2\ge Q
\end{array}\right\}.
$$
It is easily verified that $\mad_P(f,c,Q) \subset \mad_P(f),$
$\mad_L(f,c,Q) \subset \mad_L(f)$ and
$$\mad_P(f)\cap\mad_L(f)\, =   \, \bigcup_{c > 0}  (\mad_P(f,c,Q)\cap\mad_L(f,c,Q))\ . $$
For convenience we will omit the parameter $Q$ where it is
irrelevant and write $\mad_P(f,c)$ and $\mad_L(f,c)$ for
$\mad_P(f,c,Q)$ and $\mad_L(f,c,Q)$ respectively.

So it suffices to prove that the set $\mad_P(f,c)\cap\mad_L(f,c)\cap
\Ll_\theta$ has full Hausdorff dimension for some positive constant
$c$.

Geometrically, the set $\mad_P(f,c)$ consists of points that avoid
the ``neighborhood'' of each rational point $P=(p/q,r/q)$ defined by
the inequality
$$
\left|x-\frac{p}{q}\right|\left|y-\frac{r}{q}\right|<\frac{c}{f(q)q^3}.
$$
This ``neighborhood'' of $P$ will remove the interval
$\Delta(P,cf(q)^{-1})$ from $\Ll_\theta$. Without loss of generality
we can assume that $|q\theta-p|=||q\theta||$. Otherwise we just
replace the point $P$ by $P':=(p'/q,r/q)$ such that
$|q\theta-p'|=||q\theta||$. Then $\Delta(P')\supset \Delta(P)$ and
the ``neighborhood'' of~$P$ will not remove anything more than one
of $P'$.

Similarly one can show that the set $\mad_L(f,c)$ consists of points
that avoid the ``neighborhood'' of each line $L(A,B,C)$ defined by
$$
|Ax-By+C|<\frac{c}{f(|A|^*|B|^*)|A|^*|B|^*}
$$
where the coefficients $A,B,C$ satisfy $(A,B)<>(0,0)$ and
$\gcd(A,B,C)=1$. For $B=0$ it leads to the following inequality:
$$
||Ax||<\frac{c}{|A|f(|A|)}.
$$
Take $c<\inf\limits_{q\in\NN} q||q\theta||$. Then this inequality is
not true for $x=\theta$, in other words the ``neighborhood'' of the
line do not remove anything from $\Ll_\theta$. Therefore it is
sufficient to consider the lines $L(A,B,C)$ with $B\neq 0$, so the
coefficients $(A,B,C)$ will satisfy \eqref{cond_abc}. Then the
``neighborhood'' of $L(A,B,C)$ will remove the interval
$\Delta(L,cf(|A|^*|B|^*)^{-1})$ from $\Ll_\theta$.

\subsection{Cantor sets}

In the proof we will use the general Cantor framework firstly
introduced in \cite{BV_mix}. Here we reproduce the definitions and
facts which will be used in later discussion. For more details we
refer to the paper \cite{BV_mix}.

Let $I$ be a closed interval in $\RR$. Let $\RRR:=(R_n)$ with $n\in
\ZZ_{\ge 0}$ be a sequence of natural numbers and $\vr:=(r_{m,n})$
with $m,n\in \ZZ_{\ge 0},\ m\le n $ be a two parameter sequence of
non-negative real numbers.

\noindent{\bf The  construction.}  We start by subdividing the
interval $I$ into $R_0$ closed intervals $I_1$  of equal length and
denote by $\II_1$ the collection of such intervals. Thus,
$$
\#\II_1  =  R_0   \qquad {\rm and } \qquad |I_1| =  R_0^{-1}\, |{\rm
I}|  \ .
$$
Next,  we remove  at most  $r_{0,0}$ intervals $I_1$ from $\II_1$ .
Note that we do not specify which intervals should be removed but
just give an upper bound on the number of intervals to be removed.
Denote by  $\JJ_1$ the resulting collection. Thus,
\begin{equation}\label{iona1}
%\#\JJ_1  \ge   \max\{0, \#\II_1  \, - \, r_{0,0}  \} \ .
\#\JJ_1  \ge    \#\II_1   -  r_{0,0}  \, .
\end{equation}
For obvious reasons, intervals in $\JJ_1$ will be referred to as
(level one) survivors.   It will be convenient to define  $\JJ_0 :=
\{I\}$. In general, for $n \ge 0$, given  a collection $\JJ_n$ we
construct a nested collection $\JJ_{n+1}$  of closed intervals
$J_{n+1}$  using the following two operations.

{\em Splitting procedure.} We subdivide each  interval $J_n\in
\JJ_n$  into $R_n$ closed sub-intervals $I_{n+1}$ of equal length
and denote by  $\II_{n+1}$ the collection of such intervals. Thus,
    $$
    \#\II_{n+1}  =  R_n \times  \#\JJ_n   \qquad {\rm and } \qquad |I_{n+1}| = R_n^{-1}\, |J_n|    \ .
    $$

{\em Removing procedure.} For each interval $J_n\in \JJ_n$ we remove
at most $r_{n,n}$ intervals $I_{n+1} \in  \II_{n+1} $ that lie
within $J_n$.  Note that  the number of intervals $I_{n+1}$ removed
is allowed to  vary amongst  the  intervals  in  $\JJ_n$. Next, for
each interval $J_{n-1}\in \JJ_{n-1}$ we additionally remove at most
$r_{n-1,n}$ intervals $I_{n+1} \in \II_{n+1}$ that lie within $
J_{n-1}$. In general, for each interval $J_{n-k}\in \JJ_{n-k}$  $(1
\le k \le n)$ we additionally remove at most $r_{n-k,n}$ intervals
$I_{n+1} \in \II_{n+1}$ that lie within $J_{n-k}$. Then the
collection $\JJ_{n+1}$ consists of all intervals
$I_{n+1}\in\II_{n+1}$ that survive after all these removing
procedures for $k=1,2,\ldots,n$. Thus, the total number of survivors
is at most
\begin{equation}\label{iona2}
\#\JJ_{n+1}\ge R_n\#\JJ_n-\sum_{k=0}^nr_{k,n}\#\JJ_k.
\end{equation}

Finally, having constructed the nested collections $\JJ_n$ of closed
intervals   we consider the limit set
$$
 \KKK (I,\RRR,\vr) :=  \bigcap_{n=1}^\infty \bigcup_{J\in
\JJ_n} J.
$$
Any set $\KKK(I,\RRR,\vr)$ which can be achieved by the procedure
described will be referred to as a {\em $(I,\RRR,\vr)$ Cantor set.}

Of course in general it can happen that for some choice of
parameters $\RRR$ and $\vr$ and some choice of removed intervals in
removing procedure the $(I,\RRR,\vr)$ Cantor set becomes empty.
However the next result shows that with some additional conditions
on the parameters the Hausdorff dimension of this set is bounded
below.

\begin{BVThm}  Given a $({\rm
I},\RRR,\vr)$ Cantor set $\KKK (I,\RRR,\vr) $, suppose that
 $R_n\ge 4$ for all $n\in\ZZ_{\ge 0}$ and  that
\begin{equation}\label{cond_th2}
\sum_{k=0}^n \left(r_{n-k,n}\prod_{i=1}^k
\left(\frac{4}{R_{n-i}}\right)\right)\le \frac{R_n}{4}.
\end{equation}
Then
$$
\dim \KKK (I,\RRR,\vr)  \ge \liminf_{n\to \infty}(1-\log_{R_n} 2).
$$
\end{BVThm}

\noindent Here we use the convention that the product term in
\eqref{cond_th2} is one when $k=0$  and by definition  $
\log_{R_n}\!2 := \log2/ \log R_n$.  The proof of Theorem~BV4 is
presented in \cite[Theorem 4]{BV_mix}.

\subsection{Duality between points and lines}

The next two propositions show that there is a `kind' of duality
between rational points $P(p,r,q)$ and lines $L(A,B,C)$. It will
play a crucial role in our proof.
\begin{proposition}\label{prop1_p}
Let $P_1(p_1,r_1,q_1), P_2(p_2,r_2,q_2)$ be two different rational
points with $p_1/q_1\neq p_2/q_2, r_1/q_1\neq r_2/q_2$ and
$0<q_1||q_2\theta||\le q_2||q_1\theta||$. Let $L(A,B,C)$ with
$(A,B,C)$ satisfying \eqref{cond_abc} be the line passing through
$P_1,P_2$. Assume that $(P_2)_\theta\in \Delta(P_1,\delta)$. Then
\begin{equation}\label{prop1_stat_p}
(P_2)_\theta\in \Delta\left(L,
\frac{\delta^2|B|}{q_2||q_1\theta||}\cdot\frac{H(P_2)}{H(P_1)}\right)\subset
\Delta\left(L,2\delta^2\frac{H(P_2)}{H(P_1)}\right).
\end{equation}
Moreover,
\begin{equation}\label{prop1_stat2_p}
H(L)\le 4\delta H(P_1)\frac{q_2^3}{q_1^3}.
\end{equation}
\end{proposition}

\begin{proposition}\label{prop1_l}
Let $L_1(A_1,B_1,C_1), L_2(A_2,B_2,C_2)$ be two lines with integer
coefficients $(A_i,B_i,C_i)$ satisfying \eqref{cond_abc} and
$|A_2B_1|\le |A_1B_2|$. Assume that they intersect at a rational
point $P(p,r,q)$ and that $L_2\cap \Ll_\theta\in
\Delta(L_1,\delta)$. Then
\begin{equation}\label{prop1_stat_l}
L_2\cap\Ll_\theta\in \Delta\left(P,\frac{\delta^2 q}{|B_2A_1|}\cdot
\frac{H(L_2)}{H(L_1)}\right)\subset
\Delta\left(P,2\delta^2\frac{H(L_2)}{H(L_1)}\right).
\end{equation}
Moreover,
\begin{equation}\label{prop1_stat2_l}
H(P)\le 4\delta H(L_1)\frac{|B_2|^3}{|B_1|^3}.
\end{equation}
\end{proposition}

%WE NEED TO CJECK THAT |q\theta - p|=||q\theta||.!!!!!!!!!

{\bf Proof of Proposition~\ref{prop1_p}.} Since $P_1,P_2\in L$ we
have the following system of equations
$$
\left\{\begin{array}{l} Ap_1-Br_1+Cq_1=0;\\
Ap_2-Br_2+Cq_2=0;\\
A\theta-B\omega+C=0
\end{array}\right.
$$
where $\omega:=\frac{A\theta+C}{B}$. Since $p_1/q_1\neq p_2/q_2$ and
$r_1/q_1\neq r_2/q_1$ we get that the coefficients $A$ and $B$ are
nonzero. Let $A':=A/d, B':=B/d$ where $d:=(A,B)$. Then by
$(A,B,C)=1$ we get that $q_1=dq_1'$ and $q_2=dq_2'$. Then the first
two equations of the system lead to
$$
A'(p_1q_2'-p_2q_1')=B'(r_1q_2'-r_2q_1').
$$
This together with $(A',B')=1$ implies $|p_1q_2'-p_2q_1'|\ge |B'|$
and $|r_1q_2'-r_2q_1'|\ge |A'|$ or
$$
\left|\frac{p_1}{q_1}-\frac{p_2}{q_2}\right|\ge
\frac{|B|}{q_1q_2},\quad \mbox{and}\quad
\left|\frac{r_1}{q_1}-\frac{r_2}{q_2}\right|\ge \frac{|A|}{q_1q_2}.
$$
The system also gives us the following equalities
$$
|A|\left|\frac{p_1}{q_1}-\theta\right|=|B|\left|\frac{r_1}{q_1}-\omega\right|\quad\mbox{and}\quad
|A|\left|\frac{p_2}{q_2}-\theta\right|=|B|\left|\frac{r_2}{q_2}-\omega\right|.
$$
The assumption $(P_2)_\theta\in \Delta(P_1,\delta)$ is equivalent to
$$
\left|\frac{r_1}{q_1}-\frac{r_2}{q_2}\right|<\frac{\delta}{q_1^2||q_1\theta||}.
$$
Finally by the triangle inequality we find that
$$
\left|\frac{p_1}{q_1}-\frac{p_2}{q_2}\right|\le
2\max\left\{\left|\frac{p_1}{q_1}-\theta\right|,\left|\frac{p_2}{q_2}-\theta\right|\right\}=
2\max\left\{\frac{||q_1\theta||}{q_1},\frac{q_2||\theta||}{q_2}\right\}.
$$
By combining all these inequalities together we get that
$$
|B|\le q_1q_2\left|\frac{p_1}{q_1}-\frac{p_2}{q_2}\right|\le
2\max\{q_2||q_1\theta||,q_1||q_2\theta||\}=2q_2||q_1\theta||;
$$
$$
|A|\le
q_1q_2\left|\frac{r_1}{q_1}-\frac{r_2}{q_2}\right|<\frac{\delta
q_2}{q_1||q_1\theta||}.
$$
Now we are ready to calculate the bound
$$
\left|\frac{r_2}{q_2}-\omega\right|=\frac{|A|}{|B|}\frac{||q_2\theta||}{q_2}<\frac{1}{|AB|}\cdot\frac{\delta^2
q_2^2\cdot ||q_2\theta||}{q_1^2||q_1\theta||^2\cdot q_2}=
\frac{1}{|A|B^2}\cdot\frac{\delta^2 |B|\cdot
H(P_1)}{q_2||q_1\theta||\cdot H(P_2)}.
$$
Then the first inclusion in \eqref{prop1_stat_p} follows
immediately. For the second one we just use calculated estimate for
$|B|$. Also by combining the bounds for $|A|$ and $|B|$ we get
an estimate for the height $H(L)$:
$$
H(L)=|A|B^2\le \frac{4\delta q_2^3||q_1\theta||}{q_1}=4\delta
H(P_1)\frac{q_2^3}{q_1^3}.
$$

This completes the proof of Proposition~~\ref{prop1_p}. \endproof

Before we start the proof of Proposition~\ref{prop1_l} let's
establish some basic facts regarding the point of intersection of
two lines $L_1(A_1,B_1,C_1), L_2(A_2,B_2,C_2)$ with integer
coefficients $(A_i,B_i,C_i)\in\ZZ^3\backslash (\{0\}^2\times\ZZ),
(A_i,B_i,C_i)=1$; $i=1,2$. These facts will be of use in further
discussion as well. An intersection $L_1\cap L_2$ is a rational
point $P(p,r,q)$ which is the solution of the following system of
equations
$$
\left\{\begin{array}{l} A_1p-B_1r+C_1q=0;\\
A_2p-B_2r+C_2q=0
\end{array}\right.
$$
which leads to the following equalities
$$
\frac{p}{q}=\frac{B_1C_2-B_2C_1}{A_1B_2-A_2B_1}  \qquad {\rm  and }
\qquad \frac{r}{q}=\frac{A_1C_2-A_2C_1}{A_1B_2-A_2B_1} \ .
$$
Therefore we get that
\begin{equation}\label{def_prq}
|B_1C_2-B_2C_1|=dp,\ |A_1C_2-A_2C_1|=dr,\ |A_1B_2-A_2B_1|=dq.
\end{equation}
where $d:=\gcd(A_1B_2-A_2B_1, B_1C_2-B_2C_1)\in\ZZ$.

Let $i\in\{1,2\}$. It is easily verified that
$$
L_i\cap\Ll_\theta = \left(\theta,
\frac{A_i\theta+C_i}{B_i}\right)=\left(\theta,
\frac{r}{q}+\frac{A_i}{B_i}\left(\theta-\frac{p}{q}\right)\right).
$$
Therefore
$$
|L_1\cap\Ll_\theta -
L_2\cap\Ll_\theta|=\left|\frac{A_1}{B_1}-\frac{A_2}{B_2}\right|\cdot\left|\theta-\frac{p}{q}\right|=\frac{d|q\theta-p|}{|B_1B_2|}.
$$
Hence
\begin{equation}\label{ineq_qthet}
|q\theta-p|= d^{-1}|B_1B_2|\cdot|L_1\cap\Ll_\theta -
L_2\cap\Ll_\theta|
\end{equation}
and
\begin{equation}\label{ineq_qgam}
|q\omega-r|=\frac{|A_1|}{|B_1|}|q\theta-p|,\quad\text{where }\
\omega:=\frac{A_1\theta+C_1}{B_1}.
\end{equation}
\vspace{1ex}

{\bf Proof of Proposition~\ref{prop1_l}.} By \eqref{def_prq} an
upper bound for $q$ is given by
$$
q= d^{-1}|A_1B_2-A_2B_1|\le 2d^{-1}\max\{|A_1B_2|,
|A_2B_1|\}=2d^{-1}|A_1B_2|.
$$
An upper bound for $|q\theta-p|$ can be derived from
\eqref{ineq_qthet} and the assumption $L_2\cap\Ll_\theta \in
\Delta(P,\delta)$:
$$
|q\theta-p|<\frac{\delta
|B_1B_2|}{d|A_1|B_1^2}=\frac{\delta|B_2|}{d|A_1B_1|}
$$
Finally we get the required bounds
$$
|L_2\cap\Ll_\theta - P_\theta|=\frac{|A_2|}{|B_2|}\cdot
\frac{|q\theta-p|}{q}<  \frac{|A_2|\cdot \delta^2
|B_2|^2}{|B_2|\cdot d^2 |A_1B_1|^2\cdot q|q\theta-p|}\le
\frac{1}{q^2||q\theta||}\cdot \frac{\delta^2q\cdot
H(L_2)}{|B_2A_1|\cdot H(L_1)}
$$
and
$$
H(P)=q^2|q\theta-p|<
4d^{-2}|A_1B_2|^2\cdot\frac{\delta|B_2|}{d|A_1B_1|}\le 4\delta
H(L_1)\cdot\frac{|B_2|^3}{|B_1|^3}.
$$
To get the last inclusion in \eqref{prop1_stat_l} we just use
calculated bound for $q$. This completes the proof of
Proposition~\ref{prop1_l}.
\endproof

As we will see the duality between points and lines will appear
throughout the whole paper.

\section{Proof of Theorem \ref{tsc}}

\subsection{The idea}

By definition for $\theta\in\bad$ there exists a quantity
$c(\theta)>0$ such that
$$
\inf_{q\in \NN}q||q\theta||=c(\theta).
$$
In other words, for any positive integer $q$ the following
inequality is satisfied
\begin{equation}\label{ineq_cthet}
q|q\theta-p|\ge c(\theta).
\end{equation}
Let $R \ge  e^9c^{-1}(\theta)$ be an integer. Choose constants $c$
and $c_1$ sufficiently small such that they satisfy the following
inequalities
\begin{equation}\label{ineq_c1}
2^{12}c<1,\quad 2c<R^2 c_1 c(\theta), \quad c<c(\theta)
\end{equation}
and
\begin{equation}\label{ineq_c}
2^6\max\left\{\frac{c}{R^2c_1c(\theta)},2^{11}c\right\} \frac{(\log
R+2)^2R^4}{(\log2)^2}+2^{15}c_1\frac{R^3(\log R+2)}{\log 2}<1.
\end{equation}
Finally choose the parameter $Q:=c(\theta)R^2F(2)$ where
$$
F(n):=\prod_{k=1}^n k\,[\log^* k]\;\mbox{ for }n\ge 1\quad\mbox{ and
} F(n):=1\;\mbox{ for }n\le 0.
$$

The goal is to construct a $(I,\RRR,\vr)$ Cantor type set $\KK_c$
with properly chosen parameters $I,\RRR$ and $\vr$ so that $\KK_c$
is a subset of $\mad_P(f,c,Q)\cap\mad_L(f,c,Q)$. Then we use
Theorem~BV4 to estimate its Hausdorff dimension. Let $I$ be any
interval of length $c_1$ contained within the unit interval
$\{\theta\}\times [0,1]\subset \Ll_\theta$. Define $ \JJ_0:=\{I\}$.
We are going to construct, by induction on  $n$, a collection
$\JJ_n$ of closed intervals $J_n$ such that $\JJ_n $ is nested in
$\JJ_{n-1}$; that is, each interval $J_n$ in $\JJ_n$ is contained in
some interval $J_{n-1}$ in $\JJ_{n-1}$. The length of an interval
$J_n$  will be given by
$$ |J_n|  \, :=  \,   c_1    \, R^{-n}F^{-1}(n).
$$
Moreover,  each  interval $J_n $ in $\JJ_n$ will satisfy the
conditions that
\begin{equation}  \label{cond_p}
\begin{array}{cl}J_{n} \,  \cap \, \Delta(P, cf^{-1}(q)) \, =  \, \emptyset
& \forall \ \ P(p,r,q)\in\QQ^2 \ \  \mbox{with } (p,r,q)=1,\\[1ex]
&Q<H(P) < c(\theta)R^{n-1}F(n-1)
\end{array}
\end{equation}
and
\begin{equation}  \label{cond_l}
\begin{array}{cl}J_{n} \,  \cap \, \Delta(L,cf^{-1}(|A|^*|B|^*)) \, =  \, \emptyset   \; &
\forall \ L(A,B,C)\text{ with } (A,B,C)\in\ZZ^3,\ B\neq 0,
\\[1ex]
&(A,B,C)=1, \  Q<H(L) < c(\theta)R^{n-1}F(n-1) \,
\end{array}
\end{equation}
In particular, we put
$$
\KK_{c}  =   \bigcap_{n=1}^\infty \bigcup_{J\in\JJ_n}J \ .
$$
By construction, conditions \eqref{cond_p} and \eqref{cond_l} ensure
that
$$
\KK_c  \subset \mad_P(f,c)\cap\mad_L(f,c)\cap\Ll_\theta \ .
$$

The aim of the rest of the paper is to show that $\KK_c $  is in
fact a  $(I,\RRR,\vr)$ Cantor set with   $\RRR=(R_n)$ given by
\begin{equation}\label{def_rn}
R_n:=R \, (n+1) \, [\log^*\!(n+1)]
\end{equation}
and  $\vr=(r_{m,n})$  given by
%\begin{equation}\label{def_rnn}
\begin{equation}\label{def_rrnm}
r_{m,n} \, :=  \, \left\{\begin{array}{ll} 25R\log R\cdot
n^4(\log^*\!n)^4 &\mbox{
\ if }\; m=n-3\\[2ex]
0    &\mbox{ \ otherwise. }
\end{array}\right.
\end{equation}

\noindent Then Theorem~\ref{tsc} will
follow from Theorem~BV4. Indeed for $n<3$ the condition
\eqref{cond_th2} is obviously satisfied. For $n\ge 3$ and $R\ge 2^7$
we have that the
\begin{eqnarray*}
{\rm l.h.s. \ of \ } \eqref{cond_th2}  & = &  r_{n-3,n}\cdot\frac{4^3}{R_{n-1}R_{n-2}R_{n-3}}\\[2ex]
&\le  &      \frac{  4^3}{R^3} \cdot  \frac{25R\log R\cdot n^4(\log^*n)^4 }{n(n-1)(n-2)\log^*n\log^*(n-1)\log^*(n-2)} \\[2ex]
& \le &  \frac{16\cdot25\cdot 4^3\log
R}{R^3}\cdot\frac{R(n+1)[\log^*(n+1)]}{4}\le\frac{R_n }{4} \ = \
{\rm r.h.s. \ of \ } \eqref{cond_th2}  \, .
\end{eqnarray*}
Therefore Theorem~BV4 implies that
$$
\dim \KK_c   \ge \liminf_{n\to \infty} (1-\log_{R_n}\!2)=1 \,
$$
which completes the proof of Theorem~\ref{tsc}.

%
%
%Moreover, since the intervals $J_n$  are nested, in order to
%establish the nonemptiness of $ \madl_c\cap\Ll_\theta$ it suffices to
%show that each $\JJ_n$ is non-empty; i.e.
%$$\# \JJ_n   \ge 1    \qquad    \forall  \  n = 0,1, \ldots    \  . $$

\subsection{Basic construction. Splitting into collections $C_P(n,l,k)$ and $C_L(n,l,k)$}

Now we will describe the procedure of constructing the collections
$\JJ_n$. For $n=0$, we trivially have that (\ref{cond_p}),
\eqref{cond_l} are satisfied for the sole interval $I\in \JJ_0$. The
point is that by the choice of~$Q$ there are neither points nor
lines satisfying the height condition $Q<H(P),H(L)<c(\theta)$. Then
we construct $\JJ_i, i=1,2,3$ by just subdividing each $J_{i-1}$ in
$\JJ_{i-1}$ into $R\cdot i[\log^*i]$ closed intervals of equal
length. Again for the same reason the conditions (\ref{cond_p}) and
\eqref{cond_l} are satisfied for any $J_i\in\JJ_i, i=1,2,3$. Note
that
$$
\#\JJ_i=R^iF(i),\quad i=1,2,3.
$$

%For the same reason (\ref{cond_p}),\eqref{cond_l} with $n=1$ are
%trivially satisfied for any interval $J_1$  obtained by subdividing
%each $J_0$ in $\JJ_0$ into $R$ closed intervals of equal length $c_1
%R^{-1}$. Denote by $\JJ_1 $  the resulting collection of intervals
%$J_1$.  and note that
%$$
%\# \JJ_1    =  [c_1^{-1}]  \, R  \ .
%$$

In general, given $\JJ_n$ satisfying \eqref{cond_p} and
\eqref{cond_l} we wish to construct a nested collection~$\JJ_{n+1}$
of intervals $J_{n+1}$ for which (\ref{cond_p}) and \eqref{cond_l}
are satisfied with $n$ replaced by $n+1$. By definition, any
interval $J_n$ in $\JJ_n$ avoids intervals $\Delta(P,cf^{-1}(q))$
and $\Delta(L,cf^{-1}(|A|^*|B|^*))$ arising from points and lines
with height bounded above by $c(\theta)R^{n-1}F(n-1)$. Since any
`new' interval~$J_{n+1}$ is to be nested in some $J_n$, it is enough
to show that $J_{n+1}$ avoids intervals $\Delta(P,cf^{-1}(q))$ and
$\Delta(L,cf^{-1}(|A|^*|B|^*))$ arising from points and lines with
height satisfying
\begin{equation}\label{zeq2}
c(\theta)R^{n-1}F(n-1)\le H(P),H(L)<c(\theta)R^nF(n)  \   .
\end{equation}
Denote by $C_P(n)$ the collection of all rational points satisfying
this height condition.  Formally
$$
C_P(n) := \left\{P(p,r,q)\in\QQ^2  \, : \, P \ \ {\rm satsifies \
(\ref{zeq2}) \, } \right\}   \
$$
and it is precisely this  collection of rationals that comes into
play when constructing   $\JJ_{n+1}$ from~$\JJ_{n}$. By analogy for
`lines' let
%denote by $C_L(n)$ the collection of all lines with the
%same height conditions
$$
C_L(n) := \left\{L(A,B,C)  \, : \, L \ \ {\rm satsifies \
(\ref{zeq2}) \, } \right\}   \ .
$$

We now proceed with the construction. Assume that $n\ge 3$. We
subdivide each $J_n$ in $\JJ_n$  into $R_n=[R(n+1)\log^*(n+1)]$
closed intervals $I_{n+1}$ of length
$$|I_{n+1}|=c_1 R^{-n-1}F^{-1}(n+1).$$
Denote by $\II_{n+1}$ the collection of such intervals. In view of
the nested requirement, the collection $\JJ_{n+1}$ which  we are
attempting to construct will  be a  sub-collection  of  $\II_{n+1}$.
In other words, the intervals $I_{n+1}$ represent possible
candidates for $J_{n+1}$. The goal now is simple~--- it is  to
remove those `bad' intervals $I_{n+1}$ from $\II_{n+1}$  for which
\begin{equation}  \label{svt_p}
I_{n+1} \,  \cap \, \Delta(P,cf^{-1}(q)) \, \neq  \, \emptyset   \ \
\mbox{ for some \ } P(p,r,q) \in C_P(n)
\end{equation}
or
\begin{equation}  \label{svt_l}
I_{n+1} \,  \cap \, \Delta(L,cf^{-1}(|A|^*|B|^*)) \, \neq  \,
\emptyset \ \ \mbox{ for some \ } L(A,B,C) \in C_L(n)  \ .
\end{equation}
So we define
$$
\JJ_{n+1}:=\left\{J_{n+1}\in\II_{n+1}\;:\;\begin{array}{l}J_{n+1} \,
\cap \, \Delta(P,cf^{-1}(q))=\emptyset\ \mbox{ for any }P\in
C_P(n)\\[1ex]
J_{n+1} \, \cap \, \Delta(L,cf^{-1}(|A|^*|B|^*))=\emptyset\ \mbox{
for any }L\in C_L(n).
\end{array}\right\}
$$

Consider the rational point $P(p,r,q)\in C_P(n)$.  Note that since
$q^2\ge q^2 ||q\theta||=H(q)\ge cR^{n-1}F(n-1)$, we have that
\begin{equation}\label{ineq_f}
f(q)\ge\frac12\log^* (cR^{n-1}F(n-1))\log^*\frac12\log
(c(\theta)R^{n-1}F(n-1))> \frac12n(\log^*n)^2
\end{equation}
for sufficiently large $R$. We use Stirling formula to show that for
$n\ge 3$,
$$
c(\theta)R^{n-1}F(n-1)\ge c(\theta)R^{n-1}(n-1)!>(8n)^n\quad\mbox{
for }R\ge e^9c^{-1}(\theta).
$$
Therefore the left hand side of \eqref{ineq_f} is bigger than
$$
\frac12n\log(8n)\cdot \log^* (\frac12n\log (8n))>\frac12n\log^{*2}n.
$$
Note that for any line $L(A,B,C)\in C_L(n)$ we have the analogous
bound
\begin{equation}\label{ineq_fl}
f(|A|^*|B|^*)\ge \frac12n(\log^*n)^2.
\end{equation}

For $l\in \ZZ$ we split $C_P(n)$ into sub-collections
\begin{equation}\label{def_cnkl_p}
C_P(n,l):=\left\{P(p,r,q)\in C_P(n):
\begin{array}{c}c(\theta)2^lR^{n-1}F(n-1)\le
H(P)\\[1ex]
H(P)<c(\theta)2^{l+1}R^{n-1}F(n-1) \end{array}\right\}.
\end{equation}
In view of \eqref{zeq2} we have that
\begin{equation}\label{ineq_el}
2^l<R n\log^*n
\end{equation}
so
\begin{equation}\label{ineq_l}
0\le l\le [\log_2 (Rn\log^* n)] < \log_2 R + 2\log_2 n<c_3 \log^* n.
\end{equation}
where $c_3:=(\log R+2)/\log 2$ is an absolute constant independent
on $n$ and $l$.

Additionally with $k\in\ZZ$ we split the collection $C_P(n,l)$ into
sub-collections $C_P'(n,l,k)$ such that
\begin{equation}\label{def_cnk}
C_P'(n,l,k):=\left\{P(p,r,q)\in C_P(n,l)\;:\; c(\theta)2^k\le
q||q\theta||< c(\theta)2^{k+1}\right\}.
\end{equation}
Take any $P(p,r,q)\in C'_P(n,l,k)$. In view of \eqref{ineq_cthet}
the value $k$ should be nonnegative. On the other hand one can get
an upper bound for $k$ by~\eqref{zeq2}:
\begin{equation}\label{ineq_k}
0\le k\le[\log_2(R^n F(n))]<n\log_2R+n\log_2n+n\log_2\log^* n<c_3
n\log^* n,
\end{equation}
The upshot is that for fixed $n,l$ the number of classes
$C'_P(n,l,k)$ is at most $c_3n\log^* n$.

Note that within the collection $C_P'(n,l,k)$ we have very sharp
control of the height $H(P)$. Then by \eqref{def_cnkl_p} and
\eqref{def_cnk} we also have very sharp control on the value $q$ as
well, namely
\begin{equation}\label{ineq_q}
2^{l-k-1}R^{n-1}F(n-1)<q<2^{l-k+1}R^{n-1}F(n-1).
\end{equation}

Concerning the collection $C_L(n)$ we also partition it into
sub-collections. Firstly we partition it into sub-collections
$C_L(n,l)$ such that
\begin{equation}\label{def_cnkl_l}
C_L(n,l):=\left\{L\in C_L(n):
\begin{array}{c}c(\theta)2^lR^{n-1}F(n-1)\le
H(L)\\[1ex]
H(L)<c(\theta)2^{l+1}R^{n-1}F(n-1) \end{array}\right\}.
\end{equation}
Then we split $C_L(n,l)$ into sub-collections $C_L'(n,l,k)$ such
that
\begin{equation}\label{def_cnk_l}
C_L'(n,l,k):=\{L(A,B,C)\in C_L(n,l)\;:\; 2^k\le |B|<2^{k+1}\}
\end{equation}
One can check that $l$ and $k$ satisfy the same conditions
\eqref{ineq_l} and \eqref{ineq_k} as in the case of points. Note
that within each collection we have very good control of all point
and line parameters.

The procedure of removing ``bad'' intervals from $\II_{n+1}$ will be
as follows. We will firstly remove all intervals
$I_{n+1}\in\II_{n+1}$ such that there exists a point $P\in C_P(n,0)$
which satisfy $I_{n+1}\cap \Delta(P,cf^{-1}(q))\neq\emptyset$ or
there exists a line a line $L\in C_L(n,0)$ which satisfy
$I_{n+1}\cap \Delta(L,cf^{-1}(|A|^*|B|^*))\neq\emptyset$. Then we
repeat this removing procedure for collections
$$
C_P(n,1)\text{ and}\;C_L(n,1),\ldots,\;C_P(n,c_2\log^* n)\text{
and}\;C_L(n,c_2\log^* n)
$$
in exactly this order.

We will use lexicographical order for pairs in $\ZZ^2$. That is, we
say that $(a,b)\le(c,d)$ if either $a<c$ or $a=c, b\le d$. Consider
the point $P(p,r,q)\in C_P'(n,l,k)$. If there exists a pair
$(n',l')\le(n,l)$ and a line $L(A,B,C)\in C_L(n',l')$ such that
$$
H(L)<H(P)\quad\text{and}\quad\Delta(P,cf^{-1}(q))\subset\Delta(L,cf^{-1}(|A|^*|B|^*))
$$
then such a point will not remove anything more than was removed by
a line $L$. Therefore such a point can be ignored. The same is true
if there exists a point $P'(p',r',q')\in C_P(n',l')$ such that
$$
H(P')<H(P)\quad\text{and}\quad\Delta(P,cf^{-1}(q))\subset\Delta(P',cf^{-1}(q')).
$$
Therefore instead of collection $C_P'(n,l,k)$ we can work with
$$
C_P(n,l,k):=\left\{P(p,r,q)\in C_P'(n,l,k)\;\left|\;
\begin{array}{l}
\forall (n',l')<(n,l),\\
\forall L(A,B,C)\in C_L(n',l')\;\text{with }H(L)<H(P),\\[1ex]
\forall P'(p',r',q')\in C_P(n',l')\;\text{with }H(P')<H(P)\\[1ex]
\Delta(P,cf^{-1}(q))\not\subset\Delta(L,cf^{-1}(|A|^*|B|^*)),\\[1ex]
\Delta(P,cf^{-1}(q))\not\subset\Delta(P',cf^{-1}(q')).
\end{array}\right.
\right\}
$$
By the same procedure we construct the collection $C_L(n,l,k)$ from
$C_L'(n,l,k)$. Note that by the construction of $C_P(n,l,k)$ there
exists at most one point $P(p,r,q)\in C_P(n,l,k)$ with given second
coordinate $r/q$.

\subsection{Blocks of intervals $\BB_P(J)$ and $\BB_L(J)$}

Take the maximal possible constant $c_2>0$ such that
\begin{equation}\label{ineq_c2}
c_2\le \frac{1}{2^{10}c(\theta)}\quad\text{and}\quad
\frac{R^2c_1}{c_2}\in \ZZ.
\end{equation}
Fix the triple $(n,l,k)$ and consider an arbitrary interval
$J\subset \Ll_\theta$ of length $|J|=c_22^{-l}R^{-n+1}F^{-1}(n-1)$.
Then for any $P(p,r,q)\in C_P(n,l,k)$ we have
$|\Delta(P,cf^{-1}(q))|<|J|$. Indeed this is true because
$$
|J|\ge|\Delta(P,cf^{-1}(q))| \Leftrightarrow
\frac{c_2}{2^lR^{n-1}F(n-1)}\ge \frac{2c}{f(q)H(P)}
$$
$$
\stackrel{\eqref{def_cnkl_p}}\Leftarrow \frac{c_2}{2^l
R^{n-1}F(n-1)}\ge \frac{2c}{c(\theta)2^lR^{n-1}F(n-1)\cdot f(q)}.
$$
The last inequality is true provided $c_2 c(\theta)\ge 2c$ which in
turn is true by the second inequality of \eqref{ineq_c1} and
\eqref{ineq_c2}. One can easily check that the same fact is true for
any $\Delta(L,cf^{-1}(|A|^*|B|^*))$ where $L(A,B,C)\in C_L(n,l,k)$.

\begin{lemma}\label{triang_lem_p}
Let $J$ be an interval on $\Ll_\theta$ of length $|J|=c_2 2^{-l}
R^{-n+1}F^{-1}(n-1)$. Then all rational points $P(p,r,q)\in
C_P(n,l,k)$ such that $\Delta(P,cf^{-1}(q))\cap J\neq\emptyset$ lie
on a single line.
\end{lemma}

\proof \ Consider an arbitrary point $P(p,r,q)\in C_P(n,l,k)$. Then
\begin{equation}\label{cond_pq}
\left|\theta-\frac{p}{q}\right|=\frac{H(P)}{q^3}\;\stackrel{\eqref{def_cnkl_p},\eqref{ineq_q}}<\;
\frac{c(\theta)}{2^{2l-3k-4}R^{2(n-1)}F^2(n-1)}.
\end{equation}
Suppose we have three points $P_i(p_i,r_i,q_i)\in C_P(n,l,k),
i=1,2,3$ such that $\Delta(P_i,cf^{-1}(q_i))\cap J\neq\emptyset$ and
they do not lie on a single line. Then they form a triangle which
has the area at least
$$
\area(\triangle P_1P_2P_3)\ge
\frac{1}{2q_1q_2q_3}\stackrel{\eqref{ineq_q}}\ge
\frac{1}{2^{3l-3k+4}R^{3(n-1)}F^3(n-1)}.
$$
On the other hand the first coordinates $p_i/q_i$ of the points
$P_i$ should satisfy \eqref{cond_pq} and their second coordinates
$r_i/q_i$ should lie within the interval of length
$|J|+|\Delta(P_i,cf^{-1}(q_i))|\le 2|J|$. Therefore we have the
following upper bound for the area of triangle $\triangle
P_1P_2P_3$:
$$
\area(\triangle P_1P_2P_3)< \frac{
2c_22^{-l}R^{-n+1}F^{-1}(n-1)\cdot
2c(\theta)}{2^{2l-3k-4}R^{2(n-1)}F^2(n-1)}
$$
$$
\le 2^{10}c_2c(\theta)\cdot \frac{1}{2^{3l-3k+4}R^{3(n-1)}F^3(n-1)}.
$$
Finally by \eqref{ineq_c2} we get that the last value is bounded
above by
$$
\frac{1}{2^{3l-3k+4}R^{3(n-1)}F^3(n-1)}\le \area(\triangle
P_1P_2P_3)
$$
which is impossible. So we get a contradiction. \endproof

So given interval $J$ of length $c_22^{-l}R^{n-1}F^{-1}(n-1)$ if we
have at least two points $P\in C_P(n,l,k)$ as in
Lemma~\ref{triang_lem_p} then all the points with such property will
lie on a single line $L$. We denote this line by $L_J$. If there is
at most one point $P\in C_p(n,l,k)$ as in Lemma~\ref{triang_lem_p}
then we just say that $L_J$ is undefined.

Note that $L_J$ can not be horizontal because by the construction of
$C_P(n,l,k)$ there is only one point $P(p,r,q)\in C_P(n,l,k)$ with
given second coordinate $r/q$. $L_J$ can not be vertical too.
Otherwise its equation can be written as $x=C/A,$ $\gcd(A,C)=1$.
Then by the construction of $\theta$ we have that
$$
\left|\theta-\frac{p}{q}\right|=\left|\theta-\frac{C}{A}\right|\ge
\frac{c(\theta)}{A^2}
$$
which together with \eqref{cond_pq} gives us
$$
|A|\ge 2^{l-3/2k-2}R^{n-1}F(n-1).
$$
Then by defitnition of $L_J$ there exist two points
$P_1(p_1,r_1,q_1), P_2(p_2,r_2,q_2)$ with $|r_1/q_1-r_2/q_2|<2|J|$.
However
$$
\left|\frac{r_1}{q_1}-\frac{r_2}{q_2}\right|\ge
\frac{|A|}{q_1q_2}\stackrel{\eqref{ineq_q}}\ge
2^{-l+k/2-4}R^{-n+1}F^{-1}(n-1)>2|J|.
$$
So we get a contradiction.

The statement of Lemma~\ref{triang_lem_p} can be strengthened if we
have more than two points $P\in C_P(n,l,k)$ such that
$\Delta(P,cf^{-1}(q))\cap J\neq\emptyset$.

\begin{lemma}\label{lem_mtriang_p}
Let $J$ be an interval on $\Ll_\theta$ of length $|J|=c_2 2^{-l}
R^{-n+1}F^{-1}(n-1)$. Assume that there exists a line $L_J$.
Consider the sequence of consecutive intervals
$M_i\subset\Ll_\theta$, $i\in\NN$, $|M_i|=|J|$, $M_1:=J$ and bottom
end of $M_i$ coincides with the top end of $M_{i+1}$. Define the set
$$
\PPP(J,m):=\left\{P\in C_P(n,l,k)\;:\; P\in L_J \text{ and }
\Delta(P,cf^{-1}(q))\cap \left(\bigcup_{i=1}^m
M_i\right)\neq\emptyset\right\}
$$
and the value
$$
m_P(J):=\max\{m\in \NN\;|\; \# \PPP(J,m)\ge m+1\}.
$$
Then all rational points $P\in C_P(n,l,k)$ such that
$$\Delta(P,cf^{-1}(q))\cap \left(\bigcup_{i=1}^{m_P(J)}
M_i\right)\neq\emptyset$$ lie on a line $L_J$.
\end{lemma}

{\it Remark 1.} Since the number of points $P\in C_P(n,l,k), P\in
L_J$ is finite, the value $m_P(J)$ is correctly defined. Indeed
since by assumption $\#\PPP(J,1)\ge 2$, $m+1\to\infty$ and
$\#\PPP(J,m)$ is bounded then $m(J)$ exists and is finite.

{\it Remark 2.} We define the block of intervals
$$
\BB_P(J):=\bigcup_{i=1}^{m(J)}M_i.
$$
We will work with it as with one unit. If for some interval $J$ the
line $L_J$ is undefined then we define $m(J):=1$ and $\BB_P(J):=J$.
So now $m(J)$ and $\BB_P(J)$ are well defined for all intervals~$J$
of length $c_22^{-l}R^{-n+1}F^{-1}(n-1)$.

\proof \ Is similar to the proof of Lemma~\ref{triang_lem_p}. Let
$$
\PPP(J,m(J))=(P_i(p_i,r_i,q_i))_{1\le i\le m(J)+1}
$$
where the sequence $r_i/q_i$ is ordered in ascending order. Assume
that there is a point $P(p,r,q)\in C_P(n,l,k)$ such that $P\not\in
L_J$ and $\Delta(P,cf^{-1}(q))\cap \BB_P(J)\neq\emptyset$. Then the
triangle $\Delta(PP_1P_{m(J)+1})$ is splitted into $m_P(J)$ disjoint
triangles
$$
\Delta(PP_iP_{i+1}),\quad 1\le i\le m_P(J)
$$
each of which has the area
$$
\area(\Delta(PP_iP_{i+1}))\ge \frac{1}{2qq_iq_{i+1}}.
$$
On the other hand the first coordinates of the points
$P_1,\ldots,P_{m_P(J)+1}$ and $P$ satisfy \eqref{cond_pq} and their
second coordinates lie within the interval of length at most
$(m_P(J)+1)|J|$. Therefore we have the following estimate for the
area of the triangle
$$
\frac{m_P(J)}{2^{3l-3k+4}R^{3(n-1)}F^3(n-1)}\le \area(\triangle
PP_1P_{m_P(J)+1})\le
\frac{2^{9}(m_P(J)+1)c_2c(\theta)}{2^{3l-3k+4}R^{3(n-1)}F^3(n-1)}.
$$
which is impossible since the l.h.s of this inequality is bigger
than its r.h.s.
\endproof

%We will use the following straightforward corollary of the last
%lemma
%\begin{corollary}\label{corl_triang_p}
%For any interval $J$ of length $c_12^{-l}R^{-n+1}F^{-1}(n-1)$ the
%number of points $P\in C_P(n,l,k)$ such that
%$\Delta(P,cF^{-1}(q))\cap \BB_P(J)\neq\emptyset$ is bounded above by
%$$
%m_P(J)+1\le \frac{2|\BB_P(J)|}{|J|}.
%$$
%\end{corollary}

Lemmas \ref{triang_lem_p} and \ref{lem_mtriang_p} have their full
analogues for lines $L\in C_L(n,l,k)$. However the proofs
areslightly different. We will formulate them in the next two
lemmata.

\begin{lemma}\label{triang_lem_l}
Let $J$ be an interval on $\Ll_\theta$ of length $|J|=c_2 2^{-l}
R^{-n+1}F^{-1}(n-1)$. Then all lines $L(A,B,C)\in C_L(n,l,k)$ such
that $\Delta(L,cf^{-1}(|A|^*|B|^*))\cap J\neq\emptyset$ pass through
a single rational point~$P$.
\end{lemma}

\proof \ We will use the following well-known fact. Let us have
three planar lines $L_i(A_i,B_i,C_i), i=1,2,3$ defined by equations
$A_ix-B_iy+C_i=0$. Then they intersect in one point (probably at
infinity) if and only if
$$
\det\left(\begin{array}{ccc} A_1&B_1&C_1\\
A_2&B_2&C_2\\
A_3&B_3&C_3
\end{array}\right)=0.
$$

Suppose that there are three lines $L_1,L_2,L_3\in C_L(n,l,k)$ which
do not intersect at one point but their thickenings intersect $J$.
Then
$$
\left|\det\left(\begin{array}{ccc} A_1&B_1&C_1\\
A_2&B_2&C_2\\
A_3&B_3&C_3
\end{array}\right)\right|\ge 1.
$$
On the other hand we can make a vertical shifts of $L_1,L_2,L_3$ to
the distances $\delta_i<|J|+|\Delta(L_i)|<2|J|$, $i=1,2,3$ such that
they will intersect at one point on $J$. By vertically shifting a
line to the distance $\epsilon$ we change its $C$-coordinate by the
value $B\epsilon$. Therefore we have
$$
\det\left(\begin{array}{ccc} A_1&B_1&C_1+B_1\delta_1\\
A_2&B_2&C_2+B_2\delta_2\\
A_3&B_3&C_3+B_3\delta_3
\end{array}\right)=0 \quad \Rightarrow \quad
\left|\det\left(\begin{array}{ccc} A_1&B_1&B_1\delta_1\\
A_2&B_2&B_2\delta_2\\
A_3&B_3&B_3\delta_3
\end{array}\right)\right|\ge 1.
$$
However the latter determinant is bounded above by
$$
2|J|(|B_1(A_2B_3-A_3B_2)|+|B_2(A_1B_3-A_3B_1)|+|B_3(A_1B_2-A_2B_1)|)
$$$$
\stackrel{\eqref{def_cnkl_l},\eqref{def_cnk_l}}\le 2
c_2c(\theta)2^{-l}R^{-n+1}F^{-1}(n-1)\cdot 6\cdot2^{l+3}
R^{n-1}F(n-1)\stackrel{\eqref{ineq_c2}}< 1
$$
We get a contradiction.\endproof

So given interval $J$ of length $c_22^{-l}R^{n-1}F^{-1}(n-1)$ if we
have at least two lines from $C_L(n,l,k)$ as in
Lemma~\ref{triang_lem_l} then all lines with such property will
intersect at one rational point $P$. We denote this point by $P_J$.
If there is at most one line from $C_L(n,l,k)$ as in
Lemma~\ref{triang_lem_l} then we just say that $P_J$ is undefined.

The next Lemma is a ``line'' analogue of Lemma~\ref{lem_mtriang_p}.
\begin{lemma}\label{lem_mtriang_l}
Let $J$ be an interval on $\Ll_\theta$ of length $|J|=c_2 2^{-l}
R^{-n+1}F^{-1}(n-1)$. Assume that there exists a point $P_J$.
Consider the sequence of consecutive intervals
$M_i\subset\Ll_\theta$, $i\in\NN$, $|M_i|=|J|$, $M_1:=J$ and bottom
end of $M_i$ coincides with the top end of $M_{i+1}$. Define the set
$$
\LLL(J,m):=\left\{L\in C_L(n,l,k)\;:\; P_J\in L \text{ and }
\Delta(L,cf^{-1}(|A|^*|B|^*))\cap \left(\bigcup_{i=1}^m
M_i\right)\neq\emptyset\right\}
$$
and the value
$$
m_L(J):=\max\{m\in \NN\;|\; \# \LLL(J,m)\ge m+1\}.
$$
Then all lines $L\in C_L(n,l,k)$ such that
$$\Delta(L,cf^{-1}(|A|^*|B|^*))\cap \left(\bigcup_{i=1}^{m_L(J)}
M_i\right)\neq\emptyset$$ intersect at a point $P_J$.
\end{lemma}

By analogy with Remark~1 the value $m_L(J)$ is correctly defined. We
define the block of intervals
$$
\BB_L(J):=\bigcup_{i=1}^{m_L(J)}M_i.
$$
We will work with it as with one unit. As in Remark~2 if for some
interval $J$ the point~$P_J$ is not defined then we define
$m_L(J):=1$ and $\BB_L(J):=J$.

\proof \ If $m_L(J)=1$ then this is simply the statement of
Lemma~\ref{triang_lem_l}. Now assume that $m_L(J)>1$. Let
$$
\LLL(J,m_L(J))=(L_i(A_i,B_i,C_i))_{1\le i\le m_L(J)+1}.
$$
Denote by
$$\omega_i:=\frac{A_i\theta+C_i}{B_i},\quad 1\le i\le m_L(J)+1.$$
Then all the triples $(A_i,B_i,C_i)$ lie inside the figure $F$
defined by the inequalities
$$
\begin{array}{rcl}
|A_i|=\displaystyle\frac{H(L_i)}{|B_i|^2}&\stackrel{\eqref{def_cnkl_l},\eqref{def_cnk_l}}<&
c(\theta)2^{l-2k+1}R^{n-1}F(n-1),\\[2ex]
|B_i|&\stackrel{\eqref{def_cnk_l}}<&2^{k+1}\quad\text{and}\\[1ex]
|A_i\theta-B_i\omega_1+C_i|<|B_i|\cdot|\omega_1-\omega_i|&<&c_2m_L(J)2^{k+2-l}R^{-n+1}F^{-1}(n-1).
\end{array}
$$
The volume of this figure is $16c_2c(\theta)m_L(J)$ which in view of
\eqref{ineq_c2} is smaller than $\frac{1}{6}m_L(J)$. All points
$(A_i,B_i,C_i)$ together with $(0,0,0)$ lie on the plane defined by
$A_ip-B_ir+C_iq=0$. And since $\gcd(A_i,B_i,C_i)=1$ their convex
body contains at least $m_L(J)$ disjoint triangles with vertices in
points $(A_i,B_i,C_i)$ and $(0,0,0)$.

Now suppose that there is a line $L(A,B,C)\in C_L(n,l,k)$ such that
$P_J\not\in L$ and $\Delta(L,cf^{-1}(|A|^*|B|^*))\cap
\BB_L(J)\neq\emptyset$. Then $(A,B,C)\in F$ but now this point
doesn't lie on the same plane with points $(A_i,B_i,C_i)$ and
$(0,0,0)$. Then it formes at least $m_L(J)$ disjoint tetrahedrons
with them each of which has the volume at least $1/6$. Therefore the
volume of~$F$ is bounded by
$$
\frac{1}{6}m_L(J)\le \vol (F)<\frac{1}{6}m_L(J).
$$
But the last inequality is impossible. Therefore the line $L$ has to
pass through the point $P_J$.
\endproof

\subsection{Properties of blocks $\BB_P(J)$, $\BB_L(j)$ and quantities $m_P,
m_L$}\label{sec_mp}

Take an arbitrary interval $M$ of length $c_22^{-l}R^{n-1}F(n-1)$
and consider the collection $\PPP_M$ of points $P\in C_P(n,l,k)$
such that $\Delta(P,cF^{-1}(q))\cap \BB_P(M)\neq\emptyset$. Then one
of the following cases should be true.

{\bf Case 1P.} For any interval $J\subset \BB_P(M)$ such that
$|J|=|M|$,
$$
\#\SSS_J:=\#\{P\in\PPP_M\;|\;\Delta(P,cF^{-1}(q))\cap
J\neq\emptyset\}\le 2^2.
$$

{\bf Case 2P.} There exists $J\subset \BB_P(M), |J|=|M|$ such that
$\#\SSS_J>2^2$. Then the line $L_J$ is correctly defined and
therefore $L_M=L_M(A,B,C)$ is correctly defined as well. Let the
coefficient $B$ satisfy the condition
\begin{equation}\label{cond_ub}
|B|<\frac{c(\theta)2^{k+6}}{1/2c n(\log^*n)^2}
\end{equation}

{\bf Case 3P.} There exists $J\subset \BB_P(M), |J|=|M|$ such that
$\#\SSS_J>2^2$ and
\begin{equation}\label{cond_lb}
|B|\ge\frac{c(\theta)2^{k+6}}{1/2c n(\log^*n)^2}.
\end{equation}

Consider Cases~2P and 3P. Since for any $P\in\SSS_J$ all numbers
$P_\theta$ lie inside an interval of length at most $2|J|$ there are
at least two points $P_1(p_1,r_1,q_1)$ and $P_2(p_2,r_2,q_2)$ from
$\SSS_J$ such that
$$
\left|\frac{r_1}{q_1}-\frac{r_2}{q_2}\right|<2^{-1}|J|.
$$
 Without loss of generality assume that
$q_2||q_1\theta||>q_1||q_2\theta||$. Then
\begin{equation}\label{incl_p}
(P_2)_\theta\in \Delta(P_1,
2^{-1}|J|H(P_1))\stackrel{\eqref{def_cnkl_p}}\subset \Delta(P_1,
c(\theta)c_2)\stackrel{\eqref{ineq_c2}}\subset\Delta(P_1,2^{-10}).
\end{equation}
Since $L_M$ is neither vertical nor horizontal,
Proposition~\ref{prop1_p} is applicable for $\delta=2^{-10}$. It
states that
$$
(P_2)_\theta\in \Delta\left(L_M,
\frac{2^{-20}|B|}{q_2||q_1\theta||}\cdot\frac{H(P_2)}{H(P_1)}\right)
$$
and
$$
H(L_M)\le
2^{-8}H(P_1)\frac{q_2^3}{q_1^3}\stackrel{\eqref{ineq_q}}\le
\frac12H(P_1).
$$
It shows that $L_M$ belongs to the class which within the basic
construction had been considered before considering the points from
$\PPP_M$.

Now let's consider the Case~2P. By \eqref{cond_ub},
\eqref{def_cnkl_p}, \eqref{def_cnk} and \eqref{ineq_q} the inclusion
\eqref{incl_p} implies that
$$
(P_2)_\theta\in \Delta\left(L_M, \frac{2^{-11}}{1/2c
n(\log^*n)^2}\right)\stackrel{\eqref{ineq_fl}}\subset
\Delta\left(L_M, \frac{1}{8cf(|A|^*|B|^*)}\right).
$$
Now since for any $P(p',r',q')\in C_P(n,l,k)$ the distance
$|\theta-p'/q'|$ can differ from $|\theta-p/q|$ by factor at most 4
the same thing is true for the value $|\omega-r'/q'|$. An
implication of this is that for all $P\in\PPP_M$,
$$
P_\theta\in \Delta(L_M,1/2 cf^{-1}(|A|^*|B|^*)).
$$
Whence
$$
\bigcup_{P(p,r,q)\in\PPP_M}\Delta(P,cf^{-1}(q))\subset
\Delta(L_M,cf^{-1}(|A|^*|B|^*)).
$$
However by the construction of the collection $C_P(n,l,k)$, for all
$P\in C_P(n,l,k)$ intervals  $\Delta(P,cf^{-1}(q))$ are not
contained in any interval $\Delta$ previously considered. Therefore
since $\PPP_M\subset C_P(n,l,k)$ then the set $\PPP_M$ in case 2P
should be empty --- a contradiction. Therefore the Case 2P is
impossible.

Consider the last Case~3P. Let's order all the points in
$\PPP_M=(P_i(p_i,r_i,q_i))_{1\le i\le m_L(J)+1}$ in such a way that
the sequence $p_i/q_i$ is increasing. Then we have
$$
\left|\frac{p_1}{q_1}-\frac{p_{m_P(M)+1}}{q_{m_P(M)+1}}\right|\le
\left|\theta-\frac{p_1}{q_1}\right|+
\left|\theta-\frac{p_{m_P(M)+1}}{q_{m_P(M)+1}}\right|\stackrel{\eqref{cond_pq}}\le
\frac{c(\theta)}{2^{2l-3k-5}R^{2(n-1)}F^2(n-1)}.
$$
On the other hand the smallest possible difference between
consecutive numbers $p_i/q_i$ and $p_{i+1}/q_{i+1}$ is bounded below
by
$$
\frac{p_{i+1}}{q_{i+1}}-\frac{p_i}{q_i}\ge \frac{|B|}{q_iq_{i+1}}.
$$
and therefore
$$
\left|\frac{p_1}{q_1}-\frac{p_{m_P(M)+1}}{q_{m_P(M)+1}}\right|
\stackrel{\eqref{ineq_q}}\ge \frac{|B|m_P(M)}{2^{2l-2k+2}
R^{2(n-1)}F^2(n-1)}.
$$
By combining the last two inequalities and \eqref{cond_lb} we
finally get an estimate
\begin{equation}\label{ineq_mp}
m_P(M)\le c n(\log^* n)^2.
\end{equation}

Now for the same interval $M$ define the collection $\LLL_M$ of
lines $L(A,B,C)\in C_L(n,l,k)$ such that
$\Delta(L,cF^{-1}(|A|^*|B|^*))\cap \BB_L(M)\neq\emptyset$. Consider
three different cases which will be full analogues to cases 1P, 2P
and 3P.

{\bf Case 1L.} For any interval $J\subset \BB_L(M)$ such that
$|J|=|M|$,
$$
\#\SSS_J:=\#\{L(A,B,C)\in\LLL_M\;|\;\Delta(L,cF^{-1}(|A|^*|B|^*))\cap
J\neq\emptyset\}\le 2^2.
$$

{\bf Case 2L.} There exists $J\subset \BB_L(M)$, $|J|=|M|$ such that
$\#\SSS_J>2^2$. Then the point $P_J$ is correctly defined and
therefore $P_M=P_M(p,r,q)$ is correctly defined as well. Let the
coefficient $q$ satisfy the condition
\begin{equation}\label{cond_uq}
q<\frac{c(\theta)2^{l-k+3}R^{n-1}F(n-1)}{1/2c n(\log^*n)^2}.
\end{equation}

{\bf Case 3L.} There exists $J\subset \BB_L(M), |J|=|M|$ such that
$\#\SSS>2^2$ and
\begin{equation}\label{cond_lq}
q\ge\frac{c(\theta)2^{l-k+3}R^{n-1}F(n-1)}{1/2c n(\log^*n)^2}.
\end{equation}

Consider Cases~2L and 3L. The arguments will be essentially the same
to that about Cases 2P and 3P. So one can get that there are at
least two lines $L_1(A_1,B_1,C_1)$ and $L_2(A_2,B_2,C_2)$ from
$\SSS$ such that
$$
|L_1\cap\Ll_\theta-L_2\cap\Ll_\theta|<2^{-1}|J|.
$$
Without loss of generality suppose that $|A_2B_1|<|A_1B_2|$. Then
$$
L_2\cap\Ll_\theta\in \Delta(L_1, 2^{-10}).
$$
and Proposition~\ref{prop1_l} is applicable with $\delta=2^{-10}$.
Therefore arguments analogous to those used in cases 2P, 3P give us
$$
L_2\cap\Ll_\theta\in \left(P_M,\frac{2^{-19}q}{|B_2A_1|}\right)
$$
and
$$
H(P_M)\le \frac12H(L_1).
$$
Therefore the point $P_M$ is from the class which has already been
considered before considering lines from $\LLL_M$.

Now consider the Case~2L. Then by \eqref{cond_uq},
\eqref{def_cnkl_l} and \eqref{def_cnk_l} we have that
$$
L_2\cap\Ll_\theta\in \Delta\left(P_M, \frac{2^{-14}}{1/2c
n(\log^*n)^2}\right)\stackrel{\eqref{ineq_f}}\subset
\Delta\left(P_M, \frac{1}{32cf(q)}\right).
$$
Note that for any line $L(A,B,C)\in C_L(n,l,k)$ which go through
$P_M(p,r,q)$ the distance
$$
\left|\frac{A\theta+C}{B}-\frac{r}{q}\right|=\frac{|A|}{|B|}\left|\theta-\frac{p}{q}\right|
$$
can differ by factor at most 16 from the same distance for
line~$L_2$. Therefore for all $L\in\LLL_M$,
$$
L\cap\Ll_\theta\in \Delta(P_M,1/2 cf^{-1}(q)).
$$
Whence
$$
\bigcup_{L(A,B,C)\in\LLL_M}\Delta(L,cf^{-1}(|A|^*|B|^*))\subset
\Delta(P_M,cf^{-1}(q)).
$$
However since $\LLL_M\subset C_L(n,l,k)$ we get by the construction
of $C_L(n,l,k)$ that $\LLL_M$ has to be empty --- a
contradiction. Therefore the Case~2L is impossible.

Now consider the Case~3L. Let's order all the lines in
$\LLL_M=(L_i(A_i,B_i,C_i))_{1\le i\le m_L(J)+1}$ in such a way that
the sequence of the second coordinates of $L_i\cap\Ll_\theta$ is
increasing. Then we have
$$
|L_1\cap\Ll_\theta-L_{m_L(M)+1}\cap\Ll_\theta|\le
\left|L_1\cap\Ll_\theta-\frac{r}{q}\right|+
\left|L_{m_L(J)+1}\cap\Ll_\theta-\frac{r}{q}\right|
$$
$$
=\left(\frac{|A_1|}{|B_1|}+\frac{|A_{m_L(M)+1}|}{|B_{m_L(M)+1}|}\right)\cdot\frac{|q\theta-p|}{q}\stackrel{\eqref{def_cnkl_l},\eqref{def_cnk_l}}<c(\theta)2^{l-3k+2}R^{n-1}F(n-1)\cdot\frac{|q\theta-p|}{q}.
$$
On the other hand by \eqref{ineq_qthet} and \eqref{def_cnk_l} the
smallest difference between two consecutive $L_i\cap\Ll_\theta$ and
$L_{i+1}\cap\Ll_\theta$ is at least
$$
\frac{|q\theta-p|}{|B_iB_{i+1}|}> 2^{-2k-2}|q\theta-p|
$$
and therefore
$$
|L_1\cap\Ll_\theta-L_{m_L(M)+1}\cap\Ll_\theta|>m_L(M)2^{-2k-2}|q\theta-p|.
$$
By combining the upper and lower bounds for
$|L_1\cap\Ll_\theta-L_{m_L(M)+1}\cap\Ll_\theta|$ and \eqref{cond_lq}
we finally get an estimate
\begin{equation}\label{ineq_ml}
m_L(M)\le c n(\log^* n)^2.
\end{equation}

\subsection{Final step of the proof}

Let $n\ge 3$. Fix an interval $J_{n-3}\in\JJ_{n-3}$. We will firstly
estimate the quantity
$$
\#\{P(p,r,q)\in C_P(n,l,k)\;:\; \Delta(P,cf^{-1}(q))\cap
J_{n-3}\neq\emptyset\}.
$$
Split $J_{n-3}$ into
$$
K:= c_1/c_2\cdot2^lR^2(n-1)(n-2)[\log^*(n-1)][\log^*(n-2)]
$$
subintervals $M_1,\ldots,M_K$ of equal length
$c_22^{-l}R^{-n+1}F^{-1}(n-1)$ such that the bottom endpoint of
$M_i$ coincides with the top endpoint of $M_{i+1}$ ($1\le i\le
K-1$).

We start by constructing blocks from intervals $M_1,\ldots, M_K$.
Define $B_1:=\BB_P(M_{n_1})$, $B_2:=\BB_P(M_{n_2}),\ldots,$
$B_t:=\BB_P(M_{n_t})$ in such a way that $n_1:=1$ and the bottom
endpoint of $\BB_P(M_{n_i})$ coincides with the top endpoint of
$\BB_P(M_{n_{i+1}})$. By Lemma~\ref{lem_mtriang_p} for any $1\le
i<t$ we have
\begin{equation}\label{est_bi}
\#\{P(p,r,q)\in C_P(n,l,k)\;:\; \Delta(P,cf^{-1}(q))\cap
B_i\neq\emptyset\}\le m_P(M_{n_i})+1\le 2m_P(M_{n_i}).
\end{equation}
Now let's consider the last block $B_t$. The problem is that this
block is not necessarily included in $J_{n-3}$ so we need to treat
it independently. As it was discussed in Section~\ref{sec_mp}, we
have two possible cases. In Case~1P we have that for any $i\ge n_t$
$$
\#\{P(p,r,q)\in C_P(n,l,k)\;:\; \Delta(P,cf^{-1}(q))\cap
M_i\neq\emptyset\}\le 2^2.
$$
By combining it with \eqref{est_bi} we get that
\begin{equation}\label{est_pnlk}
\#\{P(p,r,q)\in C_P(n,l,k)\;:\; \Delta(P,cf^{-1}(q))\cap
J_{n-3}\neq\emptyset\}\le 4K
\end{equation}
In Case~3P we have
$$
\#\{P(p,r,q)\in C_P(n,l,k)\;:\; \Delta(P,cf^{-1}(q))\cap
B_t\neq\emptyset\}\stackrel{\eqref{ineq_mp}}\le
cn(\log^*n)^2+1\stackrel{\eqref{ineq_c1},\eqref{ineq_c2}}<K.
$$
By combining this estimate with \eqref{est_bi} we get that
$$
\#\{P(p,r,q)\in C_P(n,l,k)\;:\; \Delta(P,cf^{-1}(q))\cap
J_{n-3}\neq\emptyset\}\le 3K<4K.
$$
%which is weaker than the bound in Case~1P.

Now estimate the number of intervals $I_{n+1}\in\II_{n+1}$ which are
removed by $\Delta(P,cf^{-1}(q))$ where $P$ is some interval from
$C_P(n,l,k)$.
\begin{eqnarray}
\#\{I_{n+1}\in \II_{n+1} \!\!\!\!\! & \!\!\!\!\! :  \!\!\!\!\! \; &
\!\!\!\!\!  I_{n+1}\cap \Delta(P,cf^{-1}(q))\neq\emptyset \}  \ \le
\
\displaystyle\frac{|\Delta(P,cf^{-1}(q))|}{|I_{n+1}|}+2  \nonumber \\[2ex] &=&
\displaystyle\frac{2cR^{n+1}F(n+1)}{c_1f(q)H(q)}+2  \nonumber \\[2ex]  &\le &
\frac{2c R^2 n(n+1) \; [\log^*\!n] \; [\log^*\!(n+1)]}{c_1
c(\theta)f(q) 2^l}+2
\nonumber \\[2ex]
&\stackrel{\eqref{ineq_f}}{<} & \displaystyle\frac{8cR^2
(n+1)}{c_1c(\theta)2^l}+2  \, .   \label{iona}
\end{eqnarray}

\noindent The upshot of the  cardinality estimates \eqref{est_pnlk}
and \eqref{iona} is that
\begin{eqnarray*}
\#\{I_{n+1}\in \II_{n+1} \!\!\!\! &\!\!\!\! :  \!\!\!\! \;& \!\!\!\!
J_{n-3}  \cap  I_{n+1}\cap \Delta(P,cf^{-1}(q))\neq\emptyset\ \mbox{
for some } P\in C_P(n,l)\}
\nonumber \\
\stackrel{\eqref{ineq_k}}\le\!\! c_3n\log^*n\cdot\#\{I_{n+1}\in
\II_{n+1} \!\!\!\! &\!\!\!\!\! :  \!\!\!\!\!\! \;& \!\!\!\! J_{n-3}
\cap I_{n+1}\cap \Delta(P,cf^{-1}(q))\neq\emptyset\ \mbox{for some }
P\in C_P(n,l,k)\}
\nonumber \\[2ex]
&\le&
(c_3n\log^*n\cdot 4K)\left(2+\frac{8cR^2(n+1)}{c_1c(\theta)2^l}\right)  \nonumber \\[2ex]
&\le  & 8R^2c_3\frac{c_1}{c_2}\cdot
2^ln^3(\log^*n)^3+2^5\frac{c_3c}{c_2c(\theta)}R^4\cdot
n^4(\log^*n)^3.
\end{eqnarray*}
By analogy we get the same estimate for
$$
\#\{I_{n+1}\in \II_{n+1} \; |  \; J_{n-3}  \cap I_{n+1}\cap
\Delta(L,cf^{-1}(|A|^*|B|^*))\neq\emptyset\ \mbox{ for some } L\in
C_L(n,l)\}.
$$
By taking lines and points together and summing over $l$ satisfying
\eqref{ineq_l} we find that
$$
\#\left\{I_{n+1}\in \II_{n+1} \; \left|  \begin{array}{l} J_{n-3}
\cap I_{n+1}\cap \Delta(P,cf^{-1}(q))\neq\emptyset\ \mbox{ for some
} P\in C_P(n)\text{ or}\\[1ex]
J_{n-3} \cap I_{n+1}\cap \Delta(L,cf^{-1}(|A|^*|B|^*))\neq\emptyset\
\mbox{ for some } L\in C_L(n)
\end{array}\right.
\right\}
$$
$$
\le 16R^2c_3\frac{c_1}{c_2}n^3(\log^*n)^3\sum_{2^l<Rn\log^*n}2^l +
2^6\frac{c_3^2c}{c_2c(\theta)}R^4n^4(\log^*n)^4.
$$
If $(2^{10}c(\theta))^{-1}>R^2c_1$ then in view of \eqref{ineq_c2}
we have that $c_2=R^2c_1$. Otherwise we have that $c_2\ge
(2^{11}c(\theta))^{-1}$ and
$$\frac{c}{c_2c(\theta)}\le 2^{11}c;\quad \frac{c_1}{c_2}\le 2^{11}
c_1.
$$
In any case the last expression is bounded by
$$
\le c_4 n^4(\log^*n)^4
$$
where
$$
c_4:=2^6\max\left\{\frac{c}{R^2c_1c(\theta)},2^{11}c\right\}
\frac{(\log
R+2)^2R^4}{(\log2)^2}+16\max\{R^{-2},2^{11}c_1\}\frac{R^3(\log
R+2)}{\log 2}
$$
(recall that $c_3=(\log R+2)/\log 2$). In view of \eqref{ineq_c} the
right hand side of this inequality is bounded by
$$25R\log R\cdot n^4(\log^* n)^4=r_{n-3,n}.$$

The upshot is that for any interval $J_{n-3}\in\JJ_{n-3}$ the number
of `bad' intervals $I_{n+1}\in\II_{n+1}$ which are to be removed is
bounded by $r_{n-3,n}$. Therefore the desired set $\KKK_c$ is indeed
a $(I,\RRR,\vr)$ Cantor type set. The proof is complete.

\section{Final remark.}

In the proof of Theorem~\ref{tsc} we showed that
$\mad_P(f,c)\cap\mad_L(f,c)\cap\Ll_\theta$ contains $(I,\RRR,\vr)$
Cantor type set. It allows us to use Theorem~5 from \cite{BV_mix}:

\begin{BVThmm}\label{th_icantor}
For each integer $1\le i \le k $, suppose we are given a Cantor set
$\KKK (I,\RRR,\vr_i) $. Then
$$
\bigcap_{i=1}^{k} \KKK (I,\RRR,\vr_i)   %\, = \, \KKK (I,\RRR,\vr)
$$ %
is a $(I,\RRR,\vr) $ Cantor set where
%$
%\vr:=(r_{n,m}) $ with $$ r_{n,m} :=\sum_{i=1}^k
%r^{(i)}_{n,m}  \, .
%$$
$$
\vr:=(r_{m,n})  \quad\mbox{with } \quad r_{m,n} :=\sum_{i=1}^k
r^{(i)}_{m,n}  \, .
$$
\end{BVThmm}

Regarding the sets of the form
$\mad_P(f)\cap\mad_L(f)\cap\Ll_\theta$, Theorem~BV5 enables us to
show that for any finite family $\theta_1,\ldots,\theta_n$ of badly
approximable numbers one can find $\alpha\in\RR$ such the following
inclusion holds simultaneously for all $1\le i\le n$:
$$
(\alpha,\theta_i)\in \mad_P(f)\cap\mad_L(f).
$$
Moreover the set of such numbers $\alpha$ is of full Hausdorff
dimension. The proof is based on intersecting the corresponding
Cantor type sets $K_c(i)$ associated with each set
$\mad_P(f,c)\cap\mad_L(f,c)\cap \Ll_{\theta_i}$ for $c$ sufficiently
small and then on applying Theorem~BV4 to the intersection. We will
leave the details to the enthusiastic reader.

We also believe that the same fact will be true for countable
collection $\{\theta_i\}$ of badly approximable numbers. However it
can not be proven with existing technique.

\vspace*{8ex}

\noindent{\bf Acknowledgements.} I would like to thank Sanju Velani
for many discussions during the last years. They finally gave rise
to many ideas which made this paper possible.

\def\cprime{$'$}

\vspace{5mm}

\noindent Dzmitry A. Badziahin: Department of Mathematics,
University of York,

\vspace{0mm}

\noindent\phantom{Dzmitry A. Badziahin: }Heslington, York, YO10 5DD,
England.

%\vspace{0mm}

\noindent\phantom{Dzmitry A. Badziahin: }e-mail: db528@york.ac.uk

\end{document}